# FAST SIMULATION OF NEW COINS FROM OLD

By Şerban Nacu and Yuval Peres[1]

*University of California, Berkeley*

Let $S \subset (0,1)$. Given a known function $f : S \to (0,1)$, we consider the problem of using independent tosses of a coin with probability of heads $p$ (where $p \in S$ is unknown) to simulate a coin with probability of heads $f(p)$. We prove that if $S$ is a closed interval and $f$ is real analytic on $S$, then $f$ has a fast simulation on $S$ (the number of $p$-coin tosses needed has exponential tails). Conversely, if a function $f$ has a fast simulation on an open set, then it is real analytic on that set.

**1. Introduction.** We consider the problem of using a coin with probability of heads $p$ ($p$ unknown) to simulate a coin with probability of heads $f(p)$, where $f$ is some known function. By this we mean the following: we are allowed to toss the original $p$-coin as many times as we want. We stop at some (almost surely) finite stopping time $N$, and depending on the outcomes of the first $N$ tosses, we declare heads or tails. We want the probability of declaring a head to be exactly $f(p)$.

This problem goes back to von Neumann's 1951 article [13], where he describes an algorithm which simulates the constant function $f(p) \equiv 1/2$. It is natural to ask whether this is possible for other functions, and in 1991 Asmussen raised the question for the function $f(p) = 2p$, where it is known that $p \in (0, 1/2)$ (see [8]). The same question was raised independently but later by Propp (see [10]).

In 1994, Keane and O'Brien [8] obtained a necessary and sufficient condition for such a simulation to be possible. Consider $f : S \to [0,1]$, where $S \subset (0,1)$. Then it is possible to simulate a coin with probability of heads $f(p)$ for all $p \in S$ if and only if $f$ is constant, or $f$ is continuous and satisfies, for some $n \geq 1$,

Received September 2003; revised January 2004.
[1]Supported in part by NSF Grants DMS-01-04073 and DMS-02-44479.

(1) $$\min(f(p), 1 - f(p)) \geq \min(p, 1-p)^n \qquad \forall\, p \in S.$$

In particular, $f(p) = 2p$ cannot be simulated on $(0, 1/2)$, since the inequality (1) cannot hold for $p$ close to $1/2$. However, if we are given $\varepsilon > 0$, then an algorithm exists to simulate a $2p$-coin from tosses of a $p$-coin for $p \in (0, 1/2 - \varepsilon)$.

The methods in [8] do not provide any estimates on the number $N$ of $p$-coin tosses needed to simulate an $f(p)$-coin. The stopping time $N$ will typically be unbounded, and for fast algorithms it should have rapidly decaying tails. For example, in von Neumann's algorithm [13], the tail probabilities satisfy $\mathbf{P}_p(N > n) \leq (p^2 + (1-p)^2)^{\lfloor n/2 \rfloor}$, so they decay exponentially in $n$.

DEFINITION 1. A function $f$ has a *fast simulation* on $S$ if there exists an algorithm which simulates $f$ on $S$, and for any $p \in S$ there exist constants $C > 0, \rho < 1$ (which may depend on $p$) such that the number $N$ of required inputs satisfies $\mathbf{P}_p(N > n) \leq C\rho^n$.

REMARK. If $S$ is closed and $f$ has a fast simulation on $S$, then we can choose constants $C, \rho$ not depending on $p \in S$. See Proposition 21 for a proof.

THEOREM 1. *For any $\varepsilon > 0$, the function $f(p) = 2p$ has a fast simulation on $[0, 1/2 - \varepsilon]$.*

Building on this result, we prove:

THEOREM 2. *If $f : I \to (0,1)$ is real analytic on the closed interval $I \subset (0,1)$, then it has a fast simulation on $I$. Conversely, if a function has a fast simulation, then it is real analytic on any open subset of its domain.*

As the results stated above indicate, there is a correspondence between properties of simulation algorithms and classes of functions. Table 1 summarizes the results of [8, 10] and the present paper on this correspondence. For simplicity, in this table we restrict attention to functions $f : S \mapsto T$ where $S, T$ are closed intervals in $(0, 1)$. We do not know whether the one-sided arrows in the table can be reversed.

We prove Theorem 1 in Sections 2 and 3. In Section 2 we show that simulating $f$ is equivalent to finding sequences of certain Bernstein polynomials which approximate $f$ from above and below. If the approximations are good, then the simulations are fast. In Section 3 we use this to construct a fast simulation for the function $2p$. We can do this because the Bernstein polynomials provide exponentially convergent approximations for linear functions.

In Section 4 we prove the sufficient (constructive) part of Theorem 2. This is done in several steps. First, once we have a fast simulation for $2p$, it is easy to construct fast simulations for polynomials. Using an auxiliary geometric random variable, we also obtain fast simulations for functions which have a series expansion around the origin. This proves Theorem 2 for real analytic functions that extend to an analytic function on a disk centered at the origin. For a general real analytic function, we use Möbius maps of the form $(az+b)/(cz+d)$ to map a subset of their domain to the unit disk. Since we have fast simulations for Möbius maps, this leads to fast simulations for the original function.

In particular, Theorem 2 guarantees fast simulations for any rational function $f$, over any subset of $(0,1)$ where $\varepsilon \leq f \leq 1 - \varepsilon$. This generalizes a result from [10], where the authors prove that any rational function $f:(0,1) \to (0,1)$ has a simulation by a finite automaton, which is fast.

In Section 5 we prove the necessary part of Theorem 2, and in Section 6 we describe a very simple algorithm that gives a good approximate simulation for the function $2p$ (the error decreases exponentially in the number of steps). In Section 7 we give a simple proof of the fact that any continuous function bounded away from 0 and 1 has a simulation. Finally, in Section 8 we mention some open problems.

**2. Simulation as an approximation problem.** In this section we show that a function $f$ can be simulated if and only if it can be approximated by certain polynomials, both from below and from above, and the approximations converge to $f$. Furthermore, the speed of convergence of the approximations determines the speed of the simulation (i.e., the distribution of the number of coin tosses needed).

Let $\mathbf{P}_p$ be the law of an infinite sequence $\mathbf{X} = (X_1, X_2, \ldots)$ of i.i.d. coin tosses with probability of heads $p$. By a slight abuse of notation, we also denote by $\mathbf{P}_p$ the induced law of the first $n$ tosses $X_1, \ldots, X_n$, so for $A \subset \{0,1\}^n$, $\mathbf{P}_p(A) = \mathbf{P}_p((X_1, \ldots, X_n) \in A)$.

TABLE 1

| Simulation type | | Function class | Ref. |
|---|---|---|---|
| Terminating a.s. | $\Leftrightarrow$ | $f$ continuous | [8] |
| With finite expectation | $\Rightarrow$ | $f$ Lipshitz | Proposition 23 |
| With finite $k$th moment (and uniform tails) | $\Rightarrow$ | $f \in C^k$ | Proposition 22 |
| Fast (with exponential tails) | $\Leftrightarrow$ | $f$ real analytic | Theorem 2 |
| Via pushdown automaton | $\Rightarrow$ | $f$ algebraic over $\mathbf{Q}$ | [10] |
| Via finite automaton | $\Leftrightarrow$ | $f$ rational over $\mathbf{Q}$ and $f((0,1)) \subset (0,1)$ | [10] |

Fix $n$ and consider the first $n$ tosses. Either the algorithm terminates after at most $n$ inputs (and in that case, it outputs a 1 or a 0), or it needs more than $n$ inputs. Let $A_n \subset \{0,1\}^n$ be the set of inputs where the algorithm terminates and outputs 1, and let $B_n$ be the set of inputs where either the algorithm terminates and outputs 1, or needs more than $n$ inputs. Then clearly

$$\mathbf{P}_p(A_n) \leq \mathbf{P}_p(\text{algorithm outputs } 1) \leq \mathbf{P}_p(B_n).$$

The middle term is $f(p)$. Any sequence in $\{0,1\}^n$ has probability $p^k(1-p)^{n-k}$, where $k$ is the number of 1's in the sequence, so the lower and upper bounds are polynomials of the form $\sum_k c_k p^k(1-p)^{n-k}$, with $c_k$ nonnegative integers. The probability that the algorithm needs more than $n$ inputs is $\mathbf{P}_p(B_n) - \mathbf{P}_p(A_n)$, so if the polynomials are good approximations for $f$, then the number of inputs needed has small tails.

It is less obvious that a converse also holds: given a function $f$ and a sequence of approximating polynomials with certain properties, there exists an algorithm which generates $f$, so that the probabilities of $A_n$ and $B_n$ as defined above are given by the approximating polynomials. We prove this in the rest of this section.

In order to state our result in a compact form, we introduce the following.

DEFINITION 2. Let $q(x,y), r(x,y)$ be homogeneous polynomials of equal degree with real coefficients. If all coefficients of $r - q$ are nonnegative, then we write $q \preceq r$. If in addition $q \neq r$, then we write $q \prec r$.

This defines a partial order on the set of homogeneous polynomials of two variables. If $q \preceq r$, then clearly $q(x,y) \leq r(x,y)$ for all $x, y \geq 0$. The converse does not hold; for example, $xy \leq x^2 + y^2$ for all $x, y \geq 0$, but $xy \not\preceq x^2 + y^2$.

PROPOSITION 3. *If there exists an algorithm which simulates a function $f$ on a set $S \subset (0,1)$, then for all $n \geq 1$ there exist polynomials*

$$g_n(x,y) = \sum_{k=0}^{n} \binom{n}{k} a(n,k) x^k y^{n-k}, \qquad h_n(x,y) = \sum_{k=0}^{n} \binom{n}{k} b(n,k) x^k y^{n-k}$$

*with the following properties:*

(i) $0 \leq a(n,k) \leq b(n,k) \leq 1$.
(ii) $\binom{n}{k} a(n,k)$ and $\binom{n}{k} b(n,k)$ are integers.
(iii) $\lim_n g_n(p, 1-p) = f(p) = \lim_n h_n(p, 1-p)$ for all $p \in S$.
(iv) For all $m < n$, we have $(x+y)^{n-m} g_m(x,y) \preceq g_n(x,y)$ and $h_n(x,y) \preceq (x+y)^{n-m} h_m(x,y)$.

*Conversely, if there exist such polynomials* $g_n(x,y), h_n(x,y)$ *satisfying* (i)–(iv), *then there exists an algorithm which simulates* $f$ *on* $S$, *such that the number* $N$ *of inputs needed satisfies* $\mathbf{P}_p(N > n) = h_n(p, 1-p) - g_n(p, 1-p)$.

PROOF. $\Rightarrow$ Suppose an algorithm exists, consider its first $n$ inputs, and define as above $A_n \subset \{0,1\}^n$ to be the set of inputs where the algorithm outputs 1, and $B_n \subset \{0,1\}^n$ to be the set where the algorithm outputs 1 or needs more than $n$ inputs. We also partition $A_n = \bigcup A_{n,k}$ and $B_n = \bigcup B_{n,k}$ according to the number $k$ of 1's in each word. Then every element in $A_{n,k}$ or $B_{n,k}$ has probability $p^k(1-p)^{n-k}$, so if we define

$$a(n,k) = |A_{n,k}| \Big/ \binom{n}{k}, \qquad b(n,k) = |B_{n,k}| \Big/ \binom{n}{k},$$

then

$$g_n(p, 1-p) = \mathbf{P}_p(A_n), \qquad h_n(p, 1-p) = \mathbf{P}_p(B_n).$$

Conditions (i) and (ii) are clearly satisfied, and (iii) also follows easily. As discussed above, we have $g_n(p, 1-p) \leq f(p) \leq h_n(p, 1-p)$ and $\mathbf{P}_p(N > n) = h_n(p, 1-p) - g_n(p, 1-p)$; since the algorithm terminates almost surely, the difference must converge to 0. From the definition of $A_n$ and $B_n$, it is clear that $g_n(p, 1-p)$ is an increasing sequence, and $h_n(p, 1-p)$ is decreasing.

Condition (iv) must hold because of the structure of the sets $A_n$ and $B_n$. Indeed, let $m < n$ and assume $(X_1, \ldots, X_m) \in A_m$. Then $(X_1, \ldots, X_n) \in A_n$, whatever values $X_{m+1}, \ldots, X_n$ take. To make this formal, for $E \subset \{0,1\}^m$ define

$$T_{m,n}(E) = \{(X_1, \ldots, X_n) \in \{0,1\}^n : (X_1, \ldots, X_m) \in E\}.$$

That is, $T_{m,n}(E)$ is the set obtained by taking each element in $E$ and adding at the end all possible combinations of $n - m$ zeroes and ones. Partition $T_{m,n}(E) = \bigcup T_{m,n}^k(E)$, so that all words in $T_{m,n}^k(E)$ have exactly $k$ 1's. We have $T_{m,n}(A_m) \subset A_n$, so $T_{m,n}^k(A_m) \subset A_{n,k}$, so

$$|A_{n,k}| \geq |T_{m,n}^k(A_m)| = \sum_{i=0}^{k} \binom{n-m}{k-i} |A_{m,i}|,$$

which is the same as

(2) $$\binom{n}{k} a(n,k) \geq \sum_{i=0}^{k} \binom{n-m}{k-i} \binom{m}{i} a(m,i);$$

this is equivalent to $g_n(x,y) \succeq (x+y)^m g_m(x,y)$. A similar observation holds for the sets $B_n$, and this completes the proof of (iv).

⇐ Given the numbers $a(n,k), b(n,k)$ satisfying (i)–(iv), we shall define inductively sets $A_n = \bigcup A_{n,k}$, $B_n = \bigcup B_{n,k}$ with

$$A_{n,k} \subset B_{n,k}, \qquad |A_{n,k}| = \binom{n}{k} a(n,k), \qquad |B_{n,k}| = \binom{n}{k} b(n,k).$$

We also want the extra property that if $m < n$, then $T_{m,n}(A_m) \subset A_n$ and $T_{m,n}(B_m) \supset B_n$. Then we can construct an algorithm simulating $f$ as follows: at step $n$, output 1 if in $A_n$, output 0 if in $B_n^c$, continue if in $B_n - A_n$.

We define $A_{1,0} = \{0\}$ if $a(1,0) = 1$, and $\varnothing$ otherwise. We define $A_{1,1} = \{1\}$ if $a(1,1) = 1$, and $\varnothing$ otherwise. Similarly for $B_{1,0}$ and $B_{1,1}$. Since $a(1,k) \leq b(1,k)$, we have $A_{1,k} \subset B_{1,k}$ for $k = 0, 1$. Condition (iv) guarantees that if

$$|A_{m,k}| = \binom{m}{k} a(m,k) \quad \text{and} \quad |B_{m,k}| = \binom{m}{k} b(m,k)$$

for all $k$, then

$$(3) \qquad |T_{m,n}^k(A_m)| \leq \binom{n}{k} a(n,k) \leq \binom{n}{k} b(n,k) \leq |T_{m,n}^k(B_m)|.$$

Hence we can construct the sets $A_n, B_n$ from the sets $A_m, B_m$ as follows. We want to have

$$(4) \qquad T_{m,n}^k(A_m) \subset A_{n,k} \subset B_{n,k} \subset T_{m,n}^k(B_m).$$

In view of (3), this can be done by simply choosing any total ordering of the set of binary words of length $n$ with $k$ 1's. We build $A_{n,k}$ by starting with $T_{m,n}^k(A_m)$ and then adding elements of $T_{m,n}^k(B_m)$ in increasing order until we obtain the desired cardinality $\binom{n}{k} a(n,k)$. Then we add $\binom{n}{k} b(n,k) - \binom{n}{k} a(n,k)$ extra elements to obtain $B_{n,k}$. Of course, $A_n = \bigcup A_{n,k}$ and $B_n = \bigcup B_{n,k}$. It is immediate that the sets thus defined have the desired properties, so the induction step from $m$ to $n = m + 1$ works and the proof is complete. □

REMARK A. Condition (iv) in Proposition 3 implies that the sequence $(g_n(p, 1-p))_{n \geq 1}$ is increasing, and the sequence $(h_n(p, 1-p))_{n \geq 1}$ is decreasing ( just set $x = p, y = 1 - p$).

REMARK B. It is enough to define the numbers $a(n,k)$ and $b(n,k)$ when $n$ takes values along an increasing subsequence $n_i \uparrow \infty$. Indeed, assume (iv) holds for $m = n_i, n = n_{i+1}$. Then just like above, we can construct the sets $A_n, B_n$ from the sets $A_m, B_m$ so that (4) holds. Thus we can construct inductively the sets $A_{n_i}, B_{n_i}$. The algorithm is allowed to stop only at some $n_i$; if $n_i < n < n_{i+1}$, it just continues. This amounts to defining $A_n = T_{n_i,n}(A_{n_i}), B_n = T_{n_i,n}(B_{n_i})$ for $n_i < n < n_{i+1}$. In terms of the polynomials, this means

$$g_n(x,y) = (x+y)^{n-n_i} g_{n_i}(x,y), \qquad h_n(x,y) = (x+y)^{n-n_i} h_{n_i}(x,y)$$

for $n_i < n < n_{i+1}$. This is the same as

$$a(n, k) = (k/n)a(n-1, k-1) + (1 - k/n)a(n-1, k),$$
$$b(n, k) = (k/n)b(n-1, k-1) + (1 - k/n)b(n-1, k),$$

for $n_i < n < n_{i+1}$ and all $0 \le k \le n$. In the next section we will use this for the subsequence of powers of 2, $n_i = 2^i$. Note that it is enough to check (iv) for $m = n_i, n = n_{i+1}$, because then the algorithm is well defined and (iv) must hold for all $m, n$. Similarly, it is enough to check (iii) for $n = n_i$, because the sequences $(g_n(p, 1-p))_{n \ge 1}$ and $(h_n(p, 1-p))_{n \ge 1}$ are monotone.

REMARK C. Finally, condition (ii) in Proposition 3 is not essential. Indeed, suppose we find numbers $\alpha(n, k)$ and $\beta(n, k)$ satisfying all conditions in the proposition, except for (ii). Then if we define

(5) $$a(n, k) = \left\lfloor \alpha(n, k) \binom{n}{k} \right\rfloor / \binom{n}{k}, \qquad b(n, k) = \left\lceil \beta(n, k) \binom{n}{k} \right\rceil / \binom{n}{k},$$

conditions (i) and (ii) are trivially satisfied, and (iv) is satisfied because, for arbitrary $x_i$ nonnegative reals and $c_i$ nonnegative integers,

(6) $$\left\lfloor \sum c_i x_i \right\rfloor \ge \sum c_i \lfloor x_i \rfloor, \qquad \left\lceil \sum c_i x_i \right\rceil \le \sum c_i \lceil x_i \rceil.$$

Finally, (iii) still holds for $p \ne 0, 1$ because the error introduced in $g_n$ and $h_n$ is at most $\sum_{k=0}^{n} 2p^k(1-p)^{n-k}$, which is exponentially small.

**3. Simulating linear functions.** Let $\varepsilon > 0$, and let $f(p) = (2p) \wedge (1 - 2\varepsilon)$. Since we are only interested in small $\varepsilon$, we also assume $\varepsilon < 1/8$. We will use Proposition 3 to construct an algorithm which simulates $f$. As explained in Remark B of the previous section, it is enough to define $a(n, k)$ and $b(n, k)$ when $n$ is a power of 2. Then the compatibility equations in (iv) are equivalent to

(7) $$a(2n, k)\binom{2n}{k} \ge \sum_{i=0}^{k} a(n, i)\binom{n}{i}\binom{n}{k-i},$$

(8) $$b(2n, k)\binom{2n}{k} \le \sum_{i=0}^{k} b(n, i)\binom{n}{i}\binom{n}{k-i}.$$

These can be nicely expressed in terms of the hypergeometric distribution.

DEFINITION 3. We say a random variable $X$ has hypergeometric distribution $H(2n, k, n)$ if

(9) $$\mathbf{P}(X = i) = \binom{n}{i}\binom{n}{k-i} / \binom{2n}{k}.$$

We require $0 \leq k \leq 2n$. If we have an urn with $2n$ balls of which $k$ are red, and we select a sample of $n$ balls uniformly without replacement, then $X$ is the number of red balls in the sample.

In terms of the hypergeometric, the compatibility equations (7) and (8) become

$$a(2n, k) \geq \mathbf{E} a(n, X), \tag{10}$$

$$b(2n, k) \leq \mathbf{E} b(n, X). \tag{11}$$

We will need some properties of this distribution.

LEMMA 4. *If $X$ has distribution $H(2n, k, n)$, then:*

  (i) $\mathbf{E}(X/n) = k/(2n)$.
  (ii) $\mathbf{Var}(X/n) = k(2n - k)/(4(2n - 1)n^2) \leq 1/(2n)$.
  (iii) *If $a > 0$, then* $\mathbf{P}(|X/n - k/(2n)| > a) \leq 2\exp(-2a^2 n)$.

Both (i) and (ii) are standard facts; (iii) is a standard large deviation estimate. For a proof, see, for example, [7].

Finally, we need a way to find good approximations for $f$. Proposition 3(iii) suggests we can use the Bernstein polynomials. We recall their definition and main property. See [12], Chapter 1.4 for more details.

DEFINITION 4. For any function $f:[0,1] \to \mathbf{R}$ and any integer $n > 0$, the $n$th Bernstein polynomial of $f$ is $Q_n(x) = \sum_{k=0}^n f(k/n) \binom{n}{k} x^k (1-x)^{n-k}$.

PROPOSITION 5. *If $f$ is continuous, then $Q_n(x) \to f(x)$ uniformly on $[0,1]$.*

If a function is linear on some interval, the Bernstein polynomials provide a very good approximation to it; this suggests we could use them to construct a fast algorithm for functions such as $f(p) = (2p) \wedge (1 - 2\varepsilon)$. To prove that the compatibility equations (10), (11) hold, we will need the following.

LEMMA 6. *Let $X$ be hypergeometric with distribution $H(2n, k, n)$ as defined in (9), and let $f:[0,1] \to \mathbf{R}$ be any function with $|f| \leq 1$. Then:*

  (i) *If $f$ is Lipschitz, with $|f(x) - f(y)| \leq C|x - y|$, then $|\mathbf{E}f(X/n) - f(k/(2n))| \leq C/\sqrt{2n}$.*
  (ii) *If $f$ is twice differentiable, with $|f''| \leq C$, then $|\mathbf{E}f(X/n) - f(k/(2n))| \leq C/(4n)$.*
  (iii) *If $f$ is linear on a neighborhood of $k/(2n)$, so $f(t) = Ct + D$ if $|t - k/(2n)| \leq a$, then $|\mathbf{E}f(X/n) - f(k/(2n))| \leq (2|C| + 4)\exp(-2a^2 n)$.*

PROOF. If (i) holds, then we get

$$|\mathbf{E}f(X/n) - f(k/(2n))| \le \mathbf{E}|f(X/n) - f(k/(2n))|$$
$$\le C\mathbf{E}|X/n - k/(2n)|$$
$$\le C(\mathbf{E}|X/n - k/(2n)|^2)^{1/2}$$
$$= C\mathbf{Var}(X/n)^{1/2} \le C/\sqrt{2n}.$$

If (ii) holds, then Taylor's expansion for $f$ gives

$$|f(X/n) - f(k/(2n)) - (X/n - k/(2n))f'(k/(2n))|$$
$$\le (1/2)(X/n - k/(2n))^2 \sup|f''|$$

and $\mathbf{E}(X/n - k/(2n))f'(k/(2n)) = 0$, so

$$|\mathbf{E}f(X/n) - f(k/(2n))|$$
$$= |\mathbf{E}(f(X/n) - f(k/(2n)) - (X/n - k/(2n))f'(k/(2n)))|$$
$$\le (C/2)\mathbf{E}(X/n - k/(2n))^2$$
$$= (C/2)\mathbf{Var}(X/n) \le C/(4n).$$

If (iii) holds, then let $g(t) = f(t) - Ct - D$. We have $g = 0$ on $[k/(2n) - a, k/(2n) + a]$ and $|g(t) - g(s)| \le |f(t) - f(s)| + |C||t - s| \le 2 + |C|$ $\forall t, s \in [0, 1]$. Hence

$$|\mathbf{E}f(X/n) - f(k/(2n))|$$
$$= |\mathbf{E}g(X/n) - g(k/(2n))|$$
$$\le \mathbf{E}|g(X/n) - g(k/(2n))|$$
$$= \mathbf{E}|g(X/n) - g(k/(2n))|\mathbb{1}_{|X/n-k/2n|>a}$$
$$\le (2 + |C|)\mathbf{P}(|X/n - k/(2n)| > a)$$
$$\le 2(2 + |C|)\exp(-2a^2 n).$$

This completes the proof of the lemma. □

If we specialize the lemma to $f(p) = (2p) \wedge (1 - 2\varepsilon)$, which is Lipschitz with $C = 2$ and also piecewise linear, we obtain:

PROPOSITION 7. *Let $f(p) = (2p) \wedge (1 - 2\varepsilon)$, where $\varepsilon < 1/2$. For $X$ satisfying* (9), *we have:*

(i) $|\mathbf{E}f(X/n) - f(k/(2n))| \le \sqrt{2}/\sqrt{n}$ $\forall k, n$,
(ii) $|\mathbf{E}f(X/n) - f(k/(2n))| \le 8\exp(-2\varepsilon^2 n)$ *if* $k/(2n) \le 1/2 - 2\varepsilon$.

Now we are ready to construct the algorithm. We start by defining numbers $\alpha(n,k)$, $\beta(n,k)$ which satisfy assumptions (i), (iii) and (iv) in Proposition 3 [but not (ii)]. First we prove the compatibility equations (10) and (11):

LEMMA 8. *Define*

(12) $$\alpha(n,k) = f(k/n) = (2k/n) \wedge (1 - 2\varepsilon).$$

*Then for $X$ satisfying* (9), $\alpha(2n,k) \geq \mathbf{E}\alpha(n,X)$.

PROOF. This follows from Jensen's inequality, since $f$ is concave. □

The upper bound is more complicated. We would like $\beta(n,k)$ to be close to $\alpha(n,k)$, so that the algorithm is fast. Ideally, the difference should be exponentially small. This cannot be done over the whole interval $[0,1]$, since the Bernstein polynomials do not approximate $f$ well near $1/2 - \varepsilon$, where it is not linear. To account for this, we also need a term of order $1/\sqrt{n}$, to be added if $k/n > 1/2 - 3\varepsilon$. Finally, to control the speed of the algorithm for small $p$, we also want $\beta(n,k)$ and $\alpha(n,k)$ to be in fact equal if $k/n$ is small.

To achieve this, consider the following auxiliary functions:

$$r_1(p) = C_1(p - (1/2 - 3\varepsilon))_+, \qquad r_2(p) = C_2(p - 1/9)_+.$$

The positive constants $C_1$ and $C_2$ will be determined later. Both functions are constant, equal to zero for $p$ below a certain threshold, and increase linearly above the threshold. They are continuous and convex.

LEMMA 9. *Define*

(13) $$\beta(n,k) = f(k/n) + r_1(k/n)\sqrt{2/n} + r_2(k/n)\exp(-2\varepsilon^2 n).$$

*If $\varepsilon < 1/8$ and $X$ satisfies* (9), *then* $\beta(2n,k) \leq \mathbf{E}\beta(n,X) \ \forall k,n$.

PROOF. This amounts to proving

$$f(k/(2n)) - \mathbf{E}f(X/n)$$
$$\leq \mathbf{E}r_1(X/n)\sqrt{2/n} - r_1(k/(2n))/\sqrt{2/(2n)}$$
$$+ \mathbf{E}r_2(X/n)\exp(-2\varepsilon^2 n) - r_2(k/(2n))\exp(-4\varepsilon^2 n).$$

Since $r_1$ and $r_2$ are convex, $r_1(k/(2n)) \leq \mathbf{E}r_1(X/n)$ and $r_2(k/(2n)) \leq \mathbf{E}r_2(X/n)$, so it is enough to show

$$|f(k/(2n)) - \mathbf{E}f(X/n)|$$
$$\leq r_1(k/(2n))(1 - 1/\sqrt{2})\sqrt{2/n}$$
$$+ r_2(k/(2n))\exp(-2\varepsilon^2 n)(1 - \exp(-2\varepsilon^2 n)).$$

If $k/2n \leq 1/8$, then $X/n \leq k/n \leq 1/4 \leq 1/2 - \varepsilon$, so $f(X/n) = 2X/n$ for all values of $X$, so the left-hand side is in fact zero and the inequality holds.

If $1/8 \leq k/(2n) \leq 1/2 - 2\varepsilon$, then we use the second part of Proposition 7 (the large deviation result). Thus, it suffices to show that

$$8 \leq r_2(k/(2n))(1 - \exp(-2\varepsilon^2 n)).$$

But $r_2(k/(2n)) \geq C_2(1/8 - 1/9) = C_2/72$, so it is enough to choose

$$C_2 = 72(1 - \exp(-2\varepsilon^2))^{-1}.$$

If $k/2n > 1/2 - 2\varepsilon$, we use the first part of Proposition 7. It is enough then to show that $1 \leq r_1(k/(2n))(1 - 1/\sqrt{2})$. But $r_1(k/(2n)) \geq C_1\varepsilon$, so it is enough to choose $C_1 = \varepsilon^{-1}(1 - 1/\sqrt{2})^{-1}$. This completes the proof of the lemma. $\square$

We can now restate and prove:

THEOREM 1. *For $\varepsilon \in (0, 1/8)$, the function $f(p) = 2p \wedge (1 - 2\varepsilon)$ has a simulation on $[0, 1]$, so that the number of inputs needed, $N$, satisfies $\mathbf{P}_p(N > n) \leq C\rho^n$, for all $n \geq 1$ and $p \in [0, 1/2 - 4\varepsilon]$. The constants $C$ and $\rho$ depend on $\varepsilon$ but not on $p$, and $\rho < 1$.*

PROOF. We use Proposition 3. First we prove that for $\alpha(n, k)$ and $\beta(n, k)$ defined in (12) and (13) and

$$g_n(x, y) = \sum_{k=0}^{n} \binom{n}{k} \alpha(n, k) x^k y^{n-k}, \qquad h_n(x, y) = \sum_{k=0}^{n} \binom{n}{k} \beta(n, k) x^k y^{n-k},$$

conditions (i), (iii) and (iv) are satisfied for the subsequence $n_i = 2^i$. We have already proven (iv), and as discussed in the previous section, this implies that $g_n(p, 1-p)$ is increasing and $h_n(p, 1-p)$ is decreasing. By Proposition 5, the Bernstein polynomials $g_n(p, 1-p)$ converge to $f$. Clearly, $h_n(p, 1-p) - g_n(p, 1-p) \leq \sup_k(\beta(n, k) - \alpha(n, k)) \to 0$ as $n \to \infty$, so $h_n(p, 1-p)$ also converges to $f$ and we have proven (iii). Condition (i) clearly holds for $n$ large enough.

The remaining condition (ii) does not hold for $\alpha(n, k), \beta(n, k)$, but as discussed in the previous section, we can get around this by defining

(14) $\quad a(n, k) = \left\lfloor \alpha(n, k) \binom{n}{k} \right\rfloor / \binom{n}{k}, \qquad b(n, k) = \left\lceil \beta(n, k) \binom{n}{k} \right\rceil / \binom{n}{k}.$

Note that for $k/n < 1/9$, we have $\alpha(n, k) = \beta(n, k) = 2k/n$ so $\alpha(n, k)\binom{n}{k} = 2\binom{n-1}{k-1}$ is an integer, whence $a(n, k) = b(n, k)$.

The sequences $a(n,k)$, $b(n,k)$ satisfy conditions (i)–(iv), and the tail probabilities $\mathbf{P}_p(N > n) = h_n(p, 1-p) - g_n(p, 1-p)$ satisfy

$$\mathbf{P}_p(N > n) \leq \sum_{k=0}^{n} (\beta(n,k) - \alpha(n,k)) \binom{n}{k} p^k (1-p)^{n-k} + \sum_{k=n/9}^{n} 2p^k (1-p)^{n-k}$$
(15)
$$\leq C_1 \sqrt{2/n} \sum_{k=n/2-3\varepsilon n}^{n} \binom{n}{k} p^k (1-p)^{n-k} + C_2 e^{-2\varepsilon^2 n} + \frac{2p^{n/9}}{1-p}.$$

The second term in (15) decays exponentially, and so does the third (we can use $4 \cdot 2^{-n/9}$ as an upper bound). For the first term, ignore the square root factor and look at the sum; it is equal to $\mathbf{P}(Y/n > 1/2 - 3\varepsilon)$, where $Y$ has binomial $(n, p)$ distribution. Since $p \leq 1/2 - 4\varepsilon$, a standard large deviation estimate (see [7]) guarantees that the first term in (15) is bounded above by $\exp(-2\varepsilon^2 n)$, so it also decays exponentially in $n$.

Thus we do have $\mathbf{P}_p(N > n) \leq C\rho^n$ if $n$ is a power of 2. For general $n$, write $2^k \leq n < 2^{k+1}$. Then $\mathbf{P}_p(N > n) \leq \mathbf{P}_p(N > 2^k) \leq C\rho^{2^k} \leq C(\rho^{1/2})^n$. The proof is complete. □

REMARK. Most of the proof works for a general linear function $f(p) = (ap) \wedge (1 - a\varepsilon)$, for any $a > 0$. For integer $a$ the whole proof works (with different constants). If $a$ is not an integer, then the only problem comes from rounding the coefficients; the rounding error introduced is bounded by $\sum_{k=0}^{n} p^k (1-p)^{n-k}$, which still decays exponentially, but the rate of decay approaches 1 as $p$ approaches 0. In the next section we deduce a slightly weaker version of the result for general $a$ as a consequence of the case $a = 2$.

Proposition 3 and Lemma 6 can also be used to obtain simulations for more general functions. The simulations are no longer guaranteed to be fast, but we do obtain *some* bounds for the tails of $N$:

PROPOSITION 10. *Assume $f$ satisfies $\varepsilon < f < 1 - \varepsilon$ on $(0, 1)$. Then:*

(i) *If $f$ is Lipschitz, then it can be simulated with $\mathbf{P}_p(N > n) \leq D/\sqrt{n}$ for some uniform $D > 0$.*

(ii) *If $f$ is twice differentiable, then it can be simulated with $\mathbf{P}_p(N > n) \leq D/n$ for some uniform $D > 0$.*

REMARK. Neither of these conditions guarantees that $N$ has finite expectation, though we do believe that this should be possible to achieve, at least for $C^2$ functions.

PROOF OF PROPOSITION 10. As in the proof of Theorem 1, it is enough to define numbers $\alpha(n,k)$, $\beta(n,k)$ which satisfy assumptions (i), (iii) and (iv) in Proposition 3; assumption (ii) can then be achieved by rounding as described in Remark C. We set

$$\alpha(n,k) = f(k/n) - \delta_n, \qquad \beta(n,k) = f(k/n) + \delta_n,$$

with $\delta_n \to 0$. Then (i) holds as soon as $\delta_n < \varepsilon$ and (iii) holds because $g_n(p, 1-p) = Q_n(p) - \delta_n$, $h_n(p, 1-p) = Q_n(p) + \delta_n$, where $Q_n$ are the Bernstein polynomials. It remains to check (iv), and as in the proof of Theorem 1, it is enough to do it for $m, n$ powers of 2, which amounts to checking that for hypergeometric $X$ satisfying (9), we have $\alpha(2n, k) \geq \mathbf{E}\alpha(n, X)$ and $\beta(2n, k) \leq \mathbf{E}\beta(n, X)$. From Lemma 6,

$$\alpha(2n, k) - \mathbf{E}\alpha(n, X) \geq \delta_n - \delta_{2n} - C/\sqrt{2n}$$

if $f$ is Lipschitz with constant $C$, and

$$\alpha(2n, k) - \mathbf{E}\alpha(n, X) \geq \delta_n - \delta_{2n} - C/(4n)$$

if $f$ is twice differentiable and $|f''| \leq C$. The exact same inequalities hold for $\mathbf{E}\beta(n, X) - \beta(2n, k)$. Hence we can choose $\delta_n = (1 + \sqrt{2})C/\sqrt{n}$ in the Lipschitz case, and $\delta_n = C/(2n)$ in the twice differentiable case, and the proof is complete.

□

**4. Fast simulation for other functions.** We start with some facts about random variables with exponential tails.

PROPOSITION 11. *Let $X \geq 0$ be a random variable. Then the following are equivalent:*

(i) *There exist constants $C > 0$, $\rho < 1$ such that $\mathbf{P}(X > x) \leq C\rho^x \ \forall\, x > 0$.*
(ii) *$\mathbf{E}\exp(tX) < \infty$ for some $t > 0$.*

*If these hold, we say $X$ has* exponential tails.

PROOF. Straightforward. □

PROPOSITION 12. *Let $X_i \geq 0$ be i.i.d. with exponential tails, and let $N \geq 0$ be an integer-valued random variable with exponential tails. Then $Y = X_1 + \cdots + X_N$ has exponential tails.*

PROOF. Take $t > 0$ such that $\mathbf{E}\exp(tX_1) < \infty$. Then we can find $k > 0$ such that $\rho = \mathbf{E}\exp(t(X_1 - k)) < 1$. Let $S_n = \sum_{i=1}^n X_i$. Then

$$\mathbf{P}(S_N > kn) \leq \mathbf{P}(N > n) + \mathbf{P}(S_n > kn).$$

The first term on the right-hand side decreases exponentially fast. To evaluate the second term, we use a standard large deviation estimate,

$$\mathbf{P}(S_n > kn) \leq \exp(-tkn)\mathbf{E}\exp(tS_n) = (\mathbf{E}\exp(t(X_1 - k)))^n = \rho^n,$$

so the second term also decreases exponentially fast and we are done. $\square$

REMARK. We do not assume that $N$ is independent from the $X_i$'s.

PROPOSITION 13. *Constant functions $f(p) = c \in [0,1]$ have a fast simulation on $(0,1)$.*

PROOF. For $f(p) = 1/2$, we can use von Neumann's trick: toss coins in pairs, until we obtain 10 or 01; in the first case output 1, otherwise output 0 (if we obtain 11 or 00, we toss again). We need $2N$ tosses, where $N$ has geometric distribution with parameter $p^2 + (1-p)^2$; this clearly has exponential tails (unless $p$ is 0 or 1).

For any other constant $c$, write it in base 2: $c = \sum_{n=1}^{\infty} c_n 2^{-n}$ with $c_n \in \{0,1\}$, generate fair coins using von Neumann's trick, and toss them until we get a 1. Output $c_M$, where $M$ is the number of fair coin tosses. This scheme generates $f(p) = c$, and requires $X_1 + \cdots + X_M$ $p$-coin tosses, where $X_i$ is the number of $p$-coin tosses needed to generate the $i$th fair coin. All $X_i$ have exponential tails and so does $M$, so Proposition 12 completes the proof. Note that the rate of decay of the tails depends on $p$ but not on $c$; this will be used below. $\square$

PROPOSITION 14. *Let $S, T \subset [0,1]$.*

*(i) If $f, g$ have fast simulations on $S$, then the product $f \cdot g$ has a fast simulation on $S$.*

*(ii) If $f$ has a fast simulation on $T$ and $g$ has a fast simulation on $S$, where $g(S) \subset T$, then $f \circ g$ has a fast simulation on $S$.*

*(iii) If $f, g$ have fast simulations on $S$ and $f + g < 1 - \varepsilon$ on $S$ for some $\varepsilon > 0$, then $f + g$ has a fast simulation on $S$.*

*(iv) If $f, g$ have fast simulations on $S$ and $f - g > \varepsilon$ on $S$ for some $\varepsilon > 0$, then $f - g$ has a fast simulation on $S$.*

PROOF. (i) Let $N_f, N_g$ be the number of inputs needed to simulate each function. We simulate $f$ and $g$ separately; if both algorithms output 1, we also output 1; otherwise, we output 0. This simulates $f \cdot g$ using $N_f + N_g$ inputs, which has exponential tails by Proposition 12.

(ii) We simulate $g$ using its algorithm, then feed the results to the algorithm for $f$. We need $X_1 + \cdots + X_{N_f}$ inputs, where $X_i$ are i.i.d. with the same distribution as $N_g$. This has exponential tails by Proposition 12.

(iii) We write $f + g = h \circ \psi$, where $h(p) = 2p$ and $\psi(p) = (f(p) + g(p))/2$. We proved in the previous section that $h$ has a fast simulation on $[0, (1 - \varepsilon)/2]$. To simulate $\psi$, we simulate $f$ and $g$ separately to obtain binary variables $B_f$ and $B_g$, then toss a fair coin; if the coin is heads, we output $B_f$, otherwise we output $B_g$. So $\psi$ can be simulated using $N_f + N_g + N$ inputs, where $N$ is the number of inputs needed to simulate a fair coin. Hence $\psi$ also has a fast simulation, so (iii) follows from (ii).

(iv) Clearly $f$ has a (fast) simulation iff $1 - f$ has one, so we can look at $1 - (f - g) = (1 - f) + g < 1 - \varepsilon$. The conclusion then follows from (iii). □

PROPOSITION 15. *If $a > 0$, $\varepsilon > 0$, the function $f$ has a fast simulation on $S$, and $af(p) < 1 - \varepsilon$ on $S$, then $a \cdot f$ has a fast simulation on $S$.*

PROOF. By Theorem 1, $2p$ has a fast simulation on $[0, 1/2 - \varepsilon)$. By the composition rule Proposition 14(ii), $2^n p$ has a fast simulation on $[0, 1/2^n - \varepsilon)$. For general $a > 0$, find $n$ with $a < 2^n$ and write $ap = 2^n(a/2^n)p$. We know multiplication by $2^n$ has a fast simulation; so does multiplication by $a/2^n$, because constants smaller than 1 have a fast simulation. Hence their composition $ap$ has a fast simulation on $[0, 1/a - \varepsilon)$. We apply the composition rule Proposition 14(ii) again to complete the proof. □

PROPOSITION 16. *Let $f(p) = \sum_{n=0}^{\infty} a_n p^n$ with $a_n \geq 0$ for all $n$. Let $t \in (0, 1]$ such that $f(t) < 1$. Then $f$ has a fast simulation on $[0, t - 2\varepsilon]$, $\forall \varepsilon > 0$.*

PROOF. Write
$$\frac{\varepsilon}{t} f(p) = \sum_{n=0}^{\infty} (a_n t^n) \left(\frac{p}{t - \varepsilon}\right)^n \left(\frac{t - \varepsilon}{t}\right)^n \frac{\varepsilon}{t}.$$

Since the terms $((t - \varepsilon)/t)^n (\varepsilon/t)$ are the probabilities of a geometric distribution, we can generate an $(\varepsilon/t)f(p)$-coin as follows. First we obtain $N$ with geometric distribution, so $\mathbf{P}_p(N = n) = ((t - \varepsilon)/t)^n (\varepsilon/t)$. Then we generate $N$ i.i.d. $p/(t - \varepsilon)$-coins (by Proposition 15, this can be done by a fast simulation), and we generate one $a_N t^N$-coin [since $f(t) < 1, a_N t^N < 1$]. Finally, we multiply the $N + 1$ outputs as in Proposition 14(i).

The number of coin tosses we need is $X + Y_1 + \cdots + Y_N + Z$, where $X$ is the number of tosses required to obtain $N$, $Y_i$ is the number of tosses required to generate the $i$th $p/(t - \varepsilon)$-coin, and $Z$ is the number of tosses required to generate one (constant) $a_N t^N$-coin. $Y_i$ have exponential tails by Proposition 15, and $Z$ has exponential tails (whose rate of decay does not depend on the value of $N$) by Proposition 13.

The way we obtain $N$ is we toss $(t - \varepsilon)/t$-coins until we obtain a zero; hence $X$ can itself be written as $X = W_1 + \cdots + W_N$, where $W_i$ is the number of tosses required to generate a constant $(t-\varepsilon)/t$-coin. Hence by Proposition 12, $(\varepsilon/t)f(p)$ has a fast simulation.

Finally, $f = (t/\varepsilon)(\varepsilon/t)f$ has a fast simulation by Proposition 15. □

PROPOSITION 17. *Let $f(p) = \sum_{n=0}^{\infty} a_n p^n$ have a series expansion with arbitrary coefficients $a_n \in \mathbf{R}$ and radius of convergence $R > 0$. Let $\varepsilon > 0$ and $S \subset (0,1)$ so that $\varepsilon < f < 1 - \varepsilon$ on $S$, and $\sup S < R$. Then $f$ has a fast simulation on $S$.*

PROOF. Separating the positive and negative coefficients, we can write $f = g - h$ where $g, h$ are analytic with radius of convergence at least $R$, and have nonnegative coefficients. They must also be bounded: $g \leq M$ and $h \leq M$, with $M = \sum_{n=0}^{\infty} |a_n|(\sup S)^n < \infty$. Then $g/(2M), h/(2M)$ must have fast simulations on $S$ by Proposition 16, so by Proposition 14, so does $2M(g/(2M) - h/(2M))$.
□

PROPOSITION 18. *If $f$, $g$ have fast simulations on $S$, are both bounded on $S$, $g > \varepsilon$ on $S$, and $f/g < 1 - \varepsilon$ on $S$ for some $\varepsilon > 0$, then $f/g$ has a fast simulation on $S$.*

PROOF. Let $M = \sup g$. Let $C \in (0,1)$ and $h(p) = C/(1-p) = \sum_{n=0}^{\infty} Cp^n$. By Proposition 16, this has a fast simulation on $(0, 1 - C - \varepsilon/(4M))$. We can replace $1 - p$ with $p$ by switching heads and tails; hence $\psi(p) = C/p$ has a fast simulation on $(C + \varepsilon/(4M), 1)$. Set $C = \varepsilon/(4M)$. Then $\psi$ has a fast simulation on $(\varepsilon/(2M), 1)$ and so does $g/(2M) \in (\varepsilon/(2M), 1)$, so $\psi \circ g = \varepsilon/(2g)$ has a fast simulation on $S$. So does the product $f \cdot (\psi \circ g) = (\varepsilon/2)(f/g)$, and by Proposition 15 so does $f/g$, since we know it is bounded above by $1 - \varepsilon$.
□

THEOREM 19. *Let $f$ be a real analytic function on a closed interval $[a,b] \subset (0,1)$, so $f$ is analytic on a domain $D$ containing $[a,b]$, and assume that $f(x) \in (0,1)$ for all $x \in [a,b]$. Then $f$ has a fast simulation on $[a,b]$.*

PROOF. If $D$ is the open disk of radius 1 centered at the origin, then $f$ has a series expansion with radius of convergence 1 and the result follows from Proposition 17. For a general $D$, the idea of the proof is to map one of its subdomains to the unit disk, using a map which has a fast simulation. See Figure 1.

Using a standard compactness argument, it is easy to show we can find a domain $E$ so that $[a,b] \subset E \subset D$ and $E$ is the intersection of two large open disks of equal radius. The centers of both disks are on the line $\mathrm{Re}(z) = (a + b)/2$, located symmetrically above and below the real axis. The boundaries of the disks intersect on the real axis at the points $a-t$ and $b+t$ for some small $t>0$. If we make the radius of the disks large enough, we may assume that the angle between the disks is $\pi/n$ for some large integer $n$.

We shall use a Möbius map of the form $(pz+q)/(rz+s)$ to map those disks into half-planes. Fix $c>0$. The map

$$(16) \qquad g_1(z) = \frac{c}{z-(a-t)} - \frac{c}{(b+t)-(a-t)}$$

maps the boundaries of the disks into lines going through the origin, so it maps $E$ to the domain between those two lines contained in the positive half-plane $\text{Re}(z) > 0$. The angle between the two lines is $\pi/n$, so the map $g_1^n$ maps $E$ to the positive half-plane.

The map $g_2(z) = 1 - 2/(1+z)$ maps the positive half-plane to the unit disk, so $g_2 \circ g_1^n$ maps $E$ to the unit disk. Hence $f \circ (g_1^n)^{-1} \circ (g_2)^{-1}$ is real analytic on the unit disk (it is easy to check that the inverses of $g_1^n$ and $g_2$ are analytic on their respective domains), so it has a fast simulation on any closed interval contained in $(0,1)$. It remains to check that $g_2 \circ g_1^n$ maps $[a,b]$ to such an interval, and that it has a fast simulation. Then it follows from Proposition 14(i) that $f$ also has a fast simulation.

For sufficiently large $c$, the function $g_1$ maps the interval $[a,b]$ to the interval $[g_1(b), g_1(a)]$ where $1 < g_1(b)$. Hence $1/g_1$ maps $[a,b]$ to some closed subinterval of $(0,1)$, and by Proposition 18 it has a fast simulation (as the ratio of two linear functions). Clearly, so does $1/g_1^n$. Finally, we can write $g_2 \circ g_1^n = g_3 \circ (1/g_1^n)$, where $g_3(z) = g_2(1/z) = 1 - (2z)/(1+z)$ also has a fast simulation, by the same Proposition 18. This completes the proof. □

**5. Necessary conditions for fast simulations.**

PROPOSITION 20. *Assume $f$ has a fast simulation on an open set $S \subset (0,1)$. Then $f$ is real analytic on $S$.*

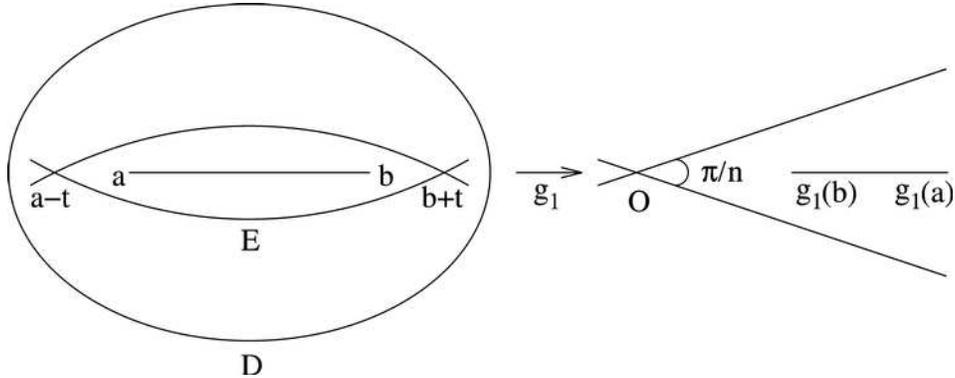

FIG. 1. *The map $g_1$.*

PROOF. Consider a fast algorithm, fix $p$ and let $f_n(p)$ be the probability that it outputs 1 after exactly $n$ steps. Then $f = \sum_{n=1}^{\infty} f_n$ and

$$0 \leq f(p) - \sum_{i=1}^{n} f_i(p) = \sum_{i=n+1}^{\infty} f_i(p) \leq C\rho^n \qquad \forall n \geq 0$$

for some constants $C > 0, \rho < 1$. Pick any $B$ with $1 < B < 1/\rho$. Since $f_n$ are polynomials, $f_n(z)$ is well defined for any complex $z$. We shall prove below that we can find $\varepsilon > 0$ so that for any complex $z$ and positive integer $n$,

(17) $$|f_n(z)| \leq B^n f_n(p) \qquad \text{if } |z - p| < \varepsilon.$$

Then for any $m > n$ and $z \in B(p, \varepsilon)$ (the open ball with center $p$ and radius $\varepsilon$), we have

$$\left| \sum_{i=n+1}^{m} f_i(z) \right| \leq \sum_{i=n+1}^{m} |f_i(z)| \leq \sum_{i=n+1}^{m} B^i f_i(p)$$

$$\leq \sum_{i=n+1}^{\infty} B^i C \rho^{i-1} = (B\rho)^n BC/(1 - B\rho).$$

Hence the sequence $\{\sum_{i=1}^{n} f_i\}$ is Cauchy on $B(p, \varepsilon)$, so it converges uniformly on $B(p, \varepsilon)$ to a limit which is analytic by a standard theorem (see [1], page 176, Theorem 1). Hence $f$ is real analytic.

To prove (17), note that $f_n$ can be written as $f_n(z) = \sum_{k=0}^{n} a_{n,k} z^k (1-z)^{n-k}$ with $a_{n,k} \geq 0$. Since $|z - p| < \varepsilon$, we have $|z| < p + \varepsilon$ and $|1 - z| < 1 - p + \varepsilon$. Choose $\varepsilon$ so $p + \varepsilon < Bp$ and $1 - p + \varepsilon < B(1 - p)$. Then

$$|z^k (1-z)^{n-k}| \leq (p + \varepsilon)^k (1 - p + \varepsilon)^{n-k} \leq B^n p^k (1-p)^{n-k}$$

and

$$\left| \sum_{k=0}^{n} a_{n,k} z^k (1-z)^{n-k} \right| \leq \sum_{k=0}^{n} a_{n,k} |z^k (1-z)^{n-k}| \leq B^n \sum_{k=0}^{n} a_{n,k} p^k (1-p)^{n-k}$$

as desired. □

PROPOSITION 21. *Assume $S \subset [0,1]$ is closed and $f$ has a fast simulation on $S$. Then the number of inputs $N$ has uniformly bounded tails: there exist constants $C, \rho$ which do not depend on $p$, so $\mathbf{P}_p(N > n) \leq C\rho^n$, $\forall p \in S$.*

PROOF. Let $g_n(p) = \mathbf{P}_p(N > n)$. Just as in Proposition 20, $g_n$ can be written as $g_n(z) = \sum_{k=0}^{n} a_{n,k} z^k (1-z)^{n-k}$ with $a_{n,k} \geq 0$, so for any $p \in (0,1)$ and $B > 1$ we can find $\varepsilon > 0$ so

(18) $$|g_n(z)| \leq B^n g_n(p) \qquad \text{if } |z - p| < \varepsilon.$$

For any $p \in S \cap (0,1)$ we have $g_n(p) \leq C_p \rho_p^n$ for some $C_p > 0, \rho_p < 1$. Setting $B = \rho_p^{-1/2}$ in (18), we obtain that there exists $\varepsilon_p > 0$ so

$$g_n(z) \leq C_p \rho_p^{n/2} \qquad \text{if } z \in (p - \varepsilon_p, p + \varepsilon_p).$$

The intervals $(p - \varepsilon_p, p + \varepsilon_p)$ cover $S$. Since $S$ is closed, it is compact, so we can find a finite subcover $(p_i - \varepsilon_{p_i}, p_i + \varepsilon_{p_i})$, $1 \leq i \leq N$. Then we can set

$$C = \max C_{p_i}, \qquad \rho = \max \rho_{p_i}^{1/2}. \qquad \square$$

REMARK. This also shows that if a function has a simulation on some $S \subset (0,1)$, then the set of $p$ where the simulation is fast is open in $S$.

PROPOSITION 22. *Assume $f$ has a simulation on an open set $S \subset (0,1)$, such that the number of inputs needed $N$ has finite $k$th moment on $S$, and furthermore the tails of the moments decrease uniformly: $\lim_{n \to \infty} \mathbf{E}_p N^k \mathbb{1}(N > n) = 0$ uniformly in $p \in S$. Then $f \in C^k(S)$ (i.e., $f$ has $k$ continuous derivatives on $S$).*

PROOF. Let $f_n$ be defined as in Proposition 20. Since $f = \sum_{n=1}^\infty f_n$, it is enough to prove that the series $\sum_{n=1}^\infty f_n^{(k)}$ converges uniformly on $S$. We shall prove that $|f_n^{(k)}| \leq Cn^k f_n$ for a uniform constant $C$. Then

$$\sum_{n=m}^\infty |f_n^{(k)}| \leq \sum_{n=m}^\infty Cn^k f_n = C\mathbf{E}_p N^k \mathbb{1}(N > m - 1)$$

converges to zero uniformly as $m \to \infty$, so the series is Cauchy and we are done. To prove the required inequality, recall that $f_n(p) = \sum_{i=0}^n a_{n,i} p^i (1 - p)^{n-i}$ with $a_{n,i} \geq 0$. Write $[i]_j = i(i-1)\cdots(i-j+1)$. From Leibniz's formula for the derivative of a product,

$$|(p^i(1-p)^{n-i})^{(k)}| = \left| \sum_{j=0}^k \binom{k}{j} (p^i)^{(j)} ((1-p)^{n-i})^{(k-j)} \right|$$

$$= \left| \sum_{j=0}^k \binom{k}{j} [i]_j p^{i-j} [n-i]_{k-j} (1-p)^{n-i-(k-j)} (-1)^{k-j} \right|$$

$$\leq \sum_{j=0}^k (k!) n^k p^i (1-p)^{n-i} / \min(p, 1-p)^k$$

$$\leq C n^k p^i (1-p)^{n-i}$$

for $C = k(k!)/\inf_{q \in B} \min(q, 1-q)^k$, where the inf is taken over some small neighborhood $B$ of $p$. It follows that $|f_n^{(k)}| \leq Cn^k f_n$ on $S$. $\square$

PROPOSITION 23. *Assume $f$ has a simulation on a closed interval $I \subset (0,1)$, such that the number of inputs needed $N$ has $\sup_{p \in I} \mathbf{E}_p(N) < \infty$. Then $f$ is Lipschitz over $I$.*

PROOF. We are given that $\mathbf{E}_p N = \sum_{n=1}^{\infty} n f_n \leq C < \infty$. Since $I$ is closed, $I \subset (\varepsilon, 1-\varepsilon)$ for some $\varepsilon$. As in the previous proposition, we obtain $|f_n'| \leq n f_n / \min(\varepsilon, 1-\varepsilon)$. Hence $|\sum_{i=1}^{n} f_i'| \leq C / \min(\varepsilon, 1-\varepsilon)$ so

$$\left| \sum_{i=1}^{n} f_i(p) - \sum_{i=1}^{n} f_i(q) \right| \leq |p-q| C / \min(\varepsilon, 1-\varepsilon).$$

Letting $n \to \infty$ completes the proof. □

**6. An approximate algorithm for doubling.** The methods described in the previous sections are essentially constructive. Proposition 3 gives a recipe for constructing an algorithm, given an approximation; all that is needed is an ordering of all binary words of length $n$ with $k$ 1's.

In the particular case of the function $f(p) = 2p$, there exists an extremely simple algorithm. It also works for any $p \in (0, 1/2)$; there is no need to bound the function away from 1. The catch is that it is approximate: it outputs 1 with probability very close to $2p$, with the error decaying exponentially in the number of steps. This must be, of course; the Keane–O'Brien results show that we could not have an *exact* algorithm with these properties. However, in practice, an approximate result may suffice.

PROPOSITION 24. *Let $p < 1/2$ and consider an asymmetric simple random walk $S_n = X_1 + \cdots + X_n$, with $\mathbf{P}_p(X_i = 1) = p = 1 - \mathbf{P}_p(X_i = -1)$. Let $A_n$ be the event that $\max(S_1, \ldots, S_n) \geq 0$. Then $\mathbf{P}_p(A_n) = \sum_{k=0}^{n} (2k/n \wedge 1) \binom{n}{k} p^k (1-p)^{n-k} = Q_n(p)$, where $Q_n$ is the $n$th Bernstein polynomial of the function $f(p) = 2p \wedge 1$.*

PROOF. We need to show that the number of paths with $k$ positive steps among the first $n$ steps, and $\max(S_1, \ldots, S_n) \geq 0$, is $(2k/n \wedge 1) \binom{n}{k}$. For $k > n/2$, this is obvious. For $k \leq n/2$, $(2k/n) \binom{n}{k} = 2 \binom{n-1}{k-1}$ and the result follows from the reflection principle (see, e.g., [3], page 197). □

Since $f$ is piecewise linear, its Bernstein polynomials converge to it exponentially fast (except at $p = 1/2$), so we obtain the following.

ALGORITHM. Run an asymmetric simple random walk $S_n = X_1 + \cdots + X_n$, with $\mathbf{P}_p(X_i = 1) = p = 1 - \mathbf{P}_p(X_i = -1)$ for at most $n$ steps. If the walk ever reaches nonnegative territory ($S_k \geq 0$ for some $1 \leq k \leq n$), output 1. Otherwise, stop after $n$ steps, output 0.

A standard large deviation estimate (see [7]) shows that if $p < 1/2$, the probability of outputting 1 is $2p - \varepsilon$, where $0 \leq \varepsilon \leq 2\exp(-2n(1/2 - p)^2)$.

See [5] for another construction of an approximate doubling algorithm.

**7. Continuous functions revisited.** In this section we use Proposition 3 to simulate any continuous function $f$ that satisfies $\varepsilon < f \leq 1 - \varepsilon$ on $(0, 1)$ for some $\varepsilon > 0$. Our proof is simpler than the original proof of [8]. However, their argument is more general since it does not assume that $f$ is bounded away from 0 and 1. We will use the following theorem of Pólya:

THEOREM 25. *Let $q(x, y)$ be a homogeneous polynomial with real coefficients satisfying $q(x, y) > 0$, $\forall x > 0, y > 0$. Then for some nonnegative integer $n$, all coefficients of $(x + y)^n q(x, y)$ are nonnegative.*

See [6], pages 57–59, for a proof. This clarifies the connection between the partial order $\preceq$ in Definition 2 and the pointwise partial order. It says that if $q(x, y) < r(x, y)$ for all $x, y > 0$, then $(x + y)^n q(x, y) \prec (x + y)^n r(x, y)$ for some $n$.

THEOREM 26 ([8]). *Let $\varepsilon > 0$ and suppose that $f : (0, 1) \mapsto [\varepsilon, 1 - \varepsilon]$ is continuous. Then $f$ admits a terminating simulation.*

PROOF. Let $i$ satisfy $2^{-i} < \varepsilon/4$. By Proposition 5, we can approximate $f - 3 \cdot 2^{-i}$ by a Bernstein polynomial $q_{m_i}$ of sufficiently high degree $m_i$ with error smaller than $2^{-i}$. More precisely,

$$q_{m_i}(x, y) = \sum_{k=0}^{m_i} \binom{m_i}{k} (f(k/m_i) - 3 \cdot 2^{-i}) x^k y^{m_i - k}$$

will satisfy $f(p) - 4 \cdot 2^{-i} < q_{m_i}(p, 1 - p) < f(p) - 2 \cdot 2^{-i}$ for all $p \in (0, 1)$.

The sequence $q_{m_i}(p, 1 - p)$ is increasing in $i$, so

$$q_{m_i}(x, y)(x + y)^{m_{i+1} - m_i} < q_{m_{i+1}}(x, y) \qquad \forall x, y > 0.$$

By Theorem 25,

$$q_{m_i}(x, y)(x + y)^{m_{i+1} - m_i + s_i} \prec q_{m_{i+1}}(x, y)(x + y)^{s_i}$$

for some integer $s_i \geq 0$. Thus if we define $n_1 = m_1$ and more generally, $n_i = m_i + (s_1 + \cdots + s_{i-1})$, then the homogeneous polynomials

$$g_{n_i}(x, y) = q_{m_i}(x, y)(x + y)^{n_i - m_i}$$

satisfy conditions (i), (iii) and (iv) in Proposition 3 along the subsequence $\{n_i\}$. Condition (ii) is easily obtained by the rounding process described in Remark C. By Remark B, once we have $g_n$ for the subsequence $n = n_i$, we can define it for all $n$. A similar construction can be used to define approximations from above $h_n$. (In fact, these approximations will require another sequence $\{s'_i\}$ analogous to $\{s_i\}$ above, and for consistency we need to use $\max\{s_i, s'_i\}$ in both approximations.) Hence by Proposition 3, $f$ has a terminating simulation algorithm. $\square$

**8. Open problems.** Theorem 2 does not settle the issue of what happens near 0 and 1, or on the boundary of the domain of analyticity of a function. An interesting example is the square root function $f(p) = \sqrt{p}$. Our methods provide fast simulations on any interval $(\varepsilon, 1]$, but if $p$ is allowed to take any value in $(0, 1)$, the best result we are aware of is the one in [10], where the authors construct a simulation using a random walk on a ladder graph. Estimates for the tails of the number of inputs needed $N$ are then given by return probabilities for a simple random walk, so $\mathbf{P}_p(N > n)$ decays like $n^{-1/2}$. We do not know whether one can do better.

QUESTION 1. Is there an algorithm that simulates $\sqrt{p}$ on $(0, 1)$, for which the number of inputs needed has finite expectation for all $p$?

REMARK. Entropy considerations (see [2], page 43) imply that if an algorithm as in Question 1 exists, then the expectation of the number of inputs cannot be uniformly bounded on $(0, 1)$. Indeed, this expectation must be at least $H(\sqrt{p})/H(p)$, where $H(p) = -p \log(p) - (1-p) \log(1-p)$ is the entropy function.

QUESTION 2. Let $J \subset (0, 1)$ be a closed interval and let $f : J \mapsto (0, 1)$ be continuous. Suppose that we have a simulation algorithm that takes as input a sequence $\{X_i\}$ of i.i.d. $p$-coins and produces a sequence of i.i.d. $f(p)$-coins. The *rate* of the algorithm (when it exists) is defined to be the limit as $n \to \infty$ of $1/n$ times the expected number of $f(p)$ coins produced from the first $n$ inputs. The rate can never exceed the entropy ratio $H(p)/H(f(p))$; see [2]. Given $J$ and $f$, are there simulation algorithms with rates arbitrarily close to the entropy ratio, uniformly for all $p \in J$?

A positive answer is known for constant $f$: for $f(p) \equiv 1/2$ variants of the von Neumann scheme (see [4, 11]) will do, and other constants follow from combining these with [9]. However, for nonconstant $f$ [except the identity and $f(p) = 1 - p$] the situation is unclear; a good example to ponder is $f(p) = p^2$.

We would also like to know whether Proposition 22 can be improved.

QUESTION 3. Is it true (possibly subject to some technical conditions) that a function has a simulation where the number of inputs has uniformly bounded $k$th moment, if and only if it has $k$ continuous derivatives?

**Acknowledgments.** We are grateful to Jim Propp for suggesting the simulation problem to us, and to Omer Angel and Elchanan Mossel for helpful discussions.

DEPARTMENT OF STATISTICS
UNIVERSITY OF CALIFORNIA
BERKELEY, CALIFORNIA 94720
USA
E-MAIL: serban@stat.berkeley.edu

DEPARTMENTS OF STATISTICS
AND MATHEMATICS
UNIVERSITY OF CALIFORNIA
BERKELEY, CALIFORNIA 94720
USA
E-MAIL: peres@stat.berkeley.edu